\documentclass[a4paper,12pt]{article}
\usepackage{amssymb,amsmath}

\newcommand{\N}{\mathbb{N}}

\newcommand{\cM}{\mathcal{M}}
\renewcommand{\b}[1]{{\bf #1}}
\newcommand{\x}{\b{x}}

\newcommand{\Z}{\mathbb{Z}}

\newcommand{\ep}{\varepsilon}
\renewcommand{\mod}{\pmod}

\newtheorem{theorem}{Theorem}
\newtheorem{lemma}{Lemma}

\begin{document}
\title{Burgess's Bounds for Character Sums}
\author{D.R. Heath-Brown\\Mathematical Institute, Oxford}
\date{}
\maketitle

% We prove that Burgess's bound gives an estimate not just for a
% single character sum but for a mean value of many such sums.  Even
% in the case of the Polya-Vinogradov inequality this seems not to
% have been previously observed.

% 11L40

\section{Introduction}

Let $\chi(n)$ be a non-principal character to modulus $q$.  Then the
well-known estimates of Burgess \cite{bur0}, \cite{bur2}, \cite{bur3}
say that if
\[S(N;H)):=\sum_{N<n\le N+H}\chi(n),\]
then for any positive integer $r\ge 2$ and any $\ep>0$ we have
\begin{equation}\label{bb}
S(N;H)\ll_{\ep,r} H^{1-1/r}q^{(r+1)/(4r^2)+\ep}
\end{equation}
uniformly in $N$, providing either that $q$ is cube free, or that $r\le 3$.
Indeed one can make the dependence on $r$
explicit, if one so wants.  Similarly the $q^{\ep}$ factor may be
replaced by a power of $d(q)\log q$ if one wishes.
The upper bound has been the best-known
for around 50 years.  The purpose of this note is to establish the following
estimate, which gives a mean-value estimate including the original
Burgess bound as a special case.

\begin{theorem}\label{mainthm}
Let $r\in\N$ and let $\ep>0$ be a real number.  Suppose that
$\chi(n)$ is a primitive character to modulus $q>1$, and let
a positive integer $H\le q$ be given.
Suppose that $0\le N_1<N_2<\ldots< N_J<q$ are integers such that
\begin{equation}\label{space}
N_{j+1}-N_j\ge H,\;\;\; (1\le j<J). 
\end{equation}
Then
\[\sum_{j=1}^J\max_{h\le
  H}|S(N_j;h)|^{3r}\ll_{\ep,r}H^{3r-3}q^{3/4+3/(4r)+\ep}\]
under any of the three conditions
\begin{enumerate}
\item[(i)] $r=1$;
\item[(ii)] $r\le 3$ and $H\ge q^{1/(2r)+\ep}$; or
\item[(iii)] $q$ is cube-free and $H\ge q^{1/(2r)+\ep}$.
\end{enumerate}
\end{theorem}

The case $J=1$ reduces to the standard Burgess estimate (which would
be trivial if one took $H\le q^{1/2r}$).  Moreover one can deduce that
there are only $O_{\ep,r}(q^{(3r+1)\ep})$ points $N_j$ for which
\[\max_{h\le  H}|S(N_j;h)|\ge H^{1-1/r}q^{(r+1)/(4r^2)-\ep},\]
for example. It would be unreasonable to ask for such a result without
the spacing condition (\ref{space}), since if $A$ and $B$ are
intervals that overlap it is possible that the behaviour of both 
$\sum_{n\in A}\chi(n)$ and $\sum_{n\in B}\chi(n)$ is affected by 
$\sum_{n\in A\cap B}\chi(n)$.

There are other results in the literature with which this estimates
should be compared.  Friedlander and Iwaniec \cite[Theorem $2'$]{FI}
establish a bound for
\[\sum_{j=1}^JS(N_j;h)\]
which can easily be used to obtain an estimate of the form
\[\sum_{j=1}^J|S(N_j;h)|^{2r}\ll_{\ep,r}h^{2r-2}q^{1/2+1/(2r)+\ep}.\]
This is superior to our result in that it involves a smaller exponent
$2r$.  However they do not include a maximum over $h$ and their 
result is subject to the condition that $h(N_J-N_1)\le q^{1+1/(2r)}$.

We should also mention the work of Chang \cite[Theorem 8]{chang}. The
result here is not so readily compared with ours, or indeed with the
Burgess estimate (\ref{bb}). However, with a certain
amount of effort one may show that our theorem  
gives a sharper bound at least when $JH^3\le q^2$.

It would have been nice to have established a result like our theorem,
but involving the $2r$-th moment.
The present methods do not allow this in
general.  However for the special case $r=1$ one can indeed achieve
this, in the following slightly more flexible form.
Specifically, suppose that $\chi(n)$ is a primitive character to 
modulus $q$, and let
$I_1,\ldots,I_J$ be disjoint subintervals of $(0,q]$.  Then for any
  $\ep>0$ we have
\begin{equation}\label{added}
\sum_{j=1}^J|\sum_{n\in I_j}\chi(n)|^2\ll_{\ep}q^{1+\ep}
\end{equation}
with an implied constant depending only on $\ep$.  This is a mild
variant of Lemma 4 of Gallagher and Montgomery \cite{GM}.
One
can deduce the P\'olya--Vinogradov as an immediate consequence of 
Lemma \ref{l3} (which is the same as Gallagher and Montgomery's Lemma
4). In fact there
are variants of (\ref{added}) for quite general character
sums. For simplicity we suppose $q$ is a prime $p$.  Let $f(x)$ and 
$g(x)$ be rational functions on $\mathbb{F}_p$, possibly identically
zero.  Then (\ref{added}) remains true if we replace $\chi(n)$
by $\chi(f(n))e_p(g(n))$, providing firstly that we exclude poles of
$f$ and $g$ from the sum, and secondly that we exclude the trivial case
in which $f$ is constant and $g$ is constant or linear.  (The implied
constant will depend on the degrees of the numerators and denominators
in $f$ and $g$.)  We leave the proof of this assertion to the reader.

For $r=1$ the ideas of this paper are closely related to those in the
article of Davenport and Erd\H{o}s \cite{DE}, which was a precursor of
Burgess's work.  For $r\ge2$ the paper follows the route to Burgess's
bounds developed in unpublished notes by Hugh Montgomery, written in
the 1970's, which were later developed into the Gallagher and Montgomery
article \cite{GM}.  In particular the mean-value lemmas in \S 2 are
essentially the same as in their paper, except that we have given the
appropriate extension to general composite moduli $q$.  
We reproduce the arguments merely for the sake of completeness.

After the mean-value lemmas in \S 2 have been established we begin the
standard attack on
the Burgess bounds in \S 3, but incorporating the sum over $N_j$ in a
non-trivial way in \S4. It is this final step that involves the 
real novelty in the paper.  This process will lead to the following key lemma.

\begin{lemma}\label{main}
Let a positive integer $r\ge 2$ and a real number $\ep>0$ be given.
Let $0\le N_1<N_2<\ldots< N_J<q$ be integers such that (\ref{space})
holds. Then for any primitive character $\chi$ to modulus $q$, and any
positive integer $H\in(q^{1/{(2r)}},q]$ we have
\begin{eqnarray}\label{est}
\lefteqn{\sum_{j=1}^J\max_{h\le H}|S(N_j;h)|^r}\nonumber\\
&&\ll_{\ep,r} 
q^{1/4+1/(4r)+\ep}H^{r-1}\{J^{2/3}+J(H^{-1}q^{1/(2r)}+Hq^{-1/2-1/(4r)})\},
\end{eqnarray}
provided either that $r\le 3$ or that $q$ is cube-free.
\end{lemma}

Throughout the paper we shall assume that $q$ is sufficiently large
in terms of $r$ and $\ep$ wherever it is convenient.  The results are
clearly trivial when $q\ll_{\ep,r}1$.  We should also point out that
we shall replace $\ep$ by a small multiple from time to time.  This
will not matter since all our results hold for all $\ep>0$.  Using
this convention we may write $q^{\ep}\log q\ll_{\ep}q^{\ep}$, for example.

\section{Preliminary Mean-value Bounds}

Our starting point, taken from previous treatments of
Burgess's bounds, is the following pair of mean value estimates.
\begin{lemma}\label{l1}
Let $r$ be a positive integer and let $\ep>0$.  Then if $\chi$ is
a primitive character to modulus $q$ we have
\[\sum_{n=1}^q|S(n;h)|^{2}\ll_{\ep} q^{1+\ep}h\]
for any $q$, while
\[\sum_{n=1}^q|S(n;h)|^{2r}\ll_{\ep,r}
  q^{\ep}(qh^r+q^{1/2}h^{2r})\]
under any of the three conditions
 \begin{enumerate}
\item[(i)] $q$ is cube-free; or 
\item[(ii)] $r=2$; or
\item[(iii)] $r=3$ and $h\le q^{1/6}$.
\end{enumerate}
\end{lemma} 
The case $r=1$ is given by Norton \cite[(2.8)]{Nort}, though the proof
is attributed to Gallagher.
For $r\ge 2$ the validity of the lemma 
under the first two conditions follows from
Burgess \cite[Lemma 8]{bur2}, using the same method as in Burgess
\cite[Lemma 8]{bur1}. The estimate under condition (iii) is given by
Burgess \cite[Theorem B]{bur3}. 

We proceed to deduce a maximal version of Lemma \ref{l1}, as in
Gallagher and Montgomery \cite[Lemma 3]{GM}.
\begin{lemma}\label{l2}
Let $r$ be a positive integer and let $\ep>0$.  Then if $\chi$ is
a primitive character to modulus $q$ and $H\in\N$ we have
\[\sum_{n=1}^q\max_{h\le H}|S(n;h)|^{2}\ll_{\ep}q^{1+\ep}H\]
for all $q$, while
\[\sum_{n=1}^q\max_{h\le H}|S(n;h)|^{2r}\ll_{\ep,r}
  q^{\ep}(qH^r+q^{1/2}H^{2r})\]
under either of the conditions
\begin{enumerate}
\item[(i)] $q$ is cube-free; or 
\item[(ii)] $2\le r\le 3$.
\end{enumerate}
\end{lemma} 

The strategy for the proof goes back to independent work of 
Rademacher \cite{Rad} and Menchov \cite{Men}, from 1922 and 1923 respectively.
It clearly suffices to consider the case in which $H=2^t$ is a power
of 2.  We will first prove the result under the assumption that $H\le
q^{1/(2r)}$. We will assume that $r\ge 2$, the case $r=1$ being similar.
Suppose that $|S(n;h)|$ attains its maximum at a positive integer
$h=h(n)\le H$, say.  We may write 
\[h=\sum_{d\in\mathcal{D}}2^{t-d}\]
for a certain set $\mathcal{D}$ of distinct non-negative integers 
$d\le t$.  Then
\[S(n;h)=\sum_{d\in\mathcal{D}}S(n+v_{n,d}2^{t-d};2^{t-d})\]
where
\[v_{n,d}=\sum_{e\in\mathcal{D},\,e<d}2^{d-e}<2^{d}.\]
By H\"{o}lder's inequality we have
\[|S(n;h)|^{2r}\le\left\{\# \mathcal{D}\right\}^{2r-1}
\left\{\sum_{d\in\mathcal{D}}|S(n+v_{n,d}2^{t-d};2^{t-d})|^{2r}\right\}.\]
We now include all possible values of $d$ and $v$ to obtain
\[|S(n;h)|^{2r}\le(t+1)^{2r-1}
\sum_{0\le d\le t}\;\sum_{0\le v<2^d}|S(n+v2^{t-d};2^{t-d})|^{2r},\]
and hence
\[\max_{h\le H}|S(n;h)|^{2r}\le(t+1)^{2r-1}
\sum_{0\le d\le t}\;\sum_{0\le v<2^d}|S(n+v2^{t-d};2^{t-d})|^{2r}.\]
We proceed to sum over $n$ modulo $q$, using Lemma \ref{l1}, and on recalling
that $H=2^t\le q^{1/(2r)}$ we deduce that
\begin{eqnarray*}
\lefteqn{\sum_{n=1}^q\max_{h\le H}|S(n;h)|^{2r}}\hspace{1cm}\\
&\ll_{\ep,r}& (t+1)^{2r-1}\sum_{0\le d\le t}\;
\sum_{0\le v<2^d}q^{\ep}(q2^{r(t-d)}+q^{1/2}2^{2r(t-d)})\\
&\ll_{\ep,r}&q^{\ep}(t+1)^{2r-1}\sum_{0\le d\le t}\;
\sum_{0\le v<2^d}(qH^r+q^{1/2}H^{2r})2^{-d}\\
&=&q^{\ep}(t+1)^{2r}(qH^r+q^{1/2}H^{2r})\\
&\ll_{\ep,r}&q^{2\ep}\rule{0cm}{0.5cm}(qH^r+q^{1/2}H^{2r}).
\end{eqnarray*}
This establishes Lemma \ref{l2} when $H$ is a power of 2 of size at
most $q^{1/(2r)}$.

To extend this to the general case, write $H_0$ for the largest power
of 2 of size at most $q^{1/(2r)}$.  Then
\[\max_{h\le H}|S(n;h)|\le
\sum_{0\le j\le H/H_0}\max_{h\le H_0}|S(n+jH_0;h)|\]
whence
\begin{eqnarray*}\sum_{n=1}^q\max_{h\le H}|S(n;h)|^{2r}&\ll& (H/H_0)^{2r-1}
\sum_{n=1}^q\;\sum_{0\le j\le H/H_0}\max_{h\le H_0}|S(n+jH_0;h)|^{2r}\\
&=&(H/H_0)^{2r-1}\sum_{0\le j\le H/H_0}\;\sum_{n=1}^q
\max_{h\le H_0}|S(n+jH_0;h)|^{2r}\\
&=&(H/H_0)^{2r-1}\sum_{0\le j\le H/H_0}\;\sum_{n\mod{q}}
\max_{h\le H_0}|S(n;h)|^{2r}\\
&\ll_{\ep,r}&(H/H_0)^{2r-1}\sum_{0\le j\le H/H_0}q^{\ep}(qH_0^r+q^{1/2}H_0^{2r})\\
&\ll_{\ep,r}&q^{\ep}(H/H_0)^{2r}(qH_0^r+q^{1/2}H_0^{2r}).
\end{eqnarray*}
However our choice of $H_0$ ensures that $qH_0^r\ll_r q^{1/2}H_0^{2r}$
and the lemma follows.

A variant of Lemma \ref{l2} allows us to sum over well spaced points.
We will only need the case $r=1$.
\begin{lemma}\label{l3}
Suppose that
$\chi(n)$ is a primitive character to modulus $q>1$, and let
a positive integer $H\le q$ be given.
Suppose that $0\le N_1<N_2<\ldots< N_J<q$ are integers satisfying the
spacing condition (\ref{space}).  Then
\[\sum_{j=1}^J\max_{h\le H}|S(N_j;h)|^{2}\ll q(\log q)^2\]
\end{lemma} 

To prove this we follow the argument in Gallagher and Montgomery
\cite[Lemma 4]{GM}.  We
first observe that for any $n\le N$ we have  
\[S(N;h)=S(n;N-n+h)-S(n;N-n).\]
If $h\le H$ it follows that
\[|S(N;h)|\le 2\max_{k\le 2H}|S(n;k)|\]
whenever $N-H<n\le N$. Then, summing over integers $n\in(N-H,N]$
we find that
\begin{equation}\label{princ}
H|S(N;h)|\le 2\sum_{n\in(N-H,N]}\max_{k\le 2H}|S(n;k)|
\end{equation}
whence H\"{o}lder's inequality yields
\[|S(N;h)|^{2r}\ll H^{-1}\sum_{n\in(N-H,N]}\max_{k\le 2H}|S(n;k)|^{2r}.\]
Since the intervals $(N_j-H,N_j]$ are disjoint modulo $q$ we then deduce 
that
\[\sum_{j=1}^J\max_{h\le H}|S(N_j;h)|^{2r}\ll H^{-1}
\sum_{n=1}^q\max_{k\le 2H}|S(n;k)|^{2r}\]
and the result follows from Lemma \ref{l2}.

We can now deduce (\ref{added}). By a dyadic subdivision it will
be enough to prove the result under the additional assumption that
there is an integer $H$ such that all the intervals $I_j$ have length
between $H/2$ and $H$. Thus we may write $I_j=(M_j,M_j+h_j]$ with
  $h_j\le H$ for $1\le j\le J$, and $M_{j+1}-M_j\ge H/2$ for $1\le
  j<J$.  We may therefore apply the case $r=1$ of Lemma \ref{l3}
  separately to the even numbered intervals and the odd numbered
  intervals to deduce (\ref{added}).

\section{Burgess's method}

In this section we will follow a mild variant of Burgess's method.
Although there are small technical differences from previous works on
the subject, there is no great novelty here.

For any prime $p<q$ which does not divide $q$ we will 
split the integers $n\in(N,N+h]$
into residue classes $n\equiv aq\mod{p}$, for $0\le a<p$.  Then we can
write $n=aq+pm$ with $m\in(N',N'+h']$ say, where
\[N'=\frac{N-aq}{p},\;\;\; h'=\frac{h}{p}.\]
We then find that
\[S(N;h)=\chi(p)\sum_{0\le a<p}S(N';h')\]
and hence
\[|S(N;h)|\le\sum_{0\le a<p}|S(N';h')|.\]
We now choose an integer parameter $P$ in the range $(\log q)^2\le P<q/2$, and 
sum the above estimate for all primes $p\in(P,2P]$ not 
dividing $q$.  Since the number
of such primes is asymptotically $P/(\log P)$ we deduce that
\begin{equation}\label{b1}
P/(\log P)|S(N;h)|\ll\sum_{P<p\le 2P}\;\sum_{0\le a<p}|S(N';h')|.
\end{equation}

We now apply the inequality (\ref{princ}), with $H$ replaced by
$H/P$. Since we have $h'\le H/P$ we deduce that
\[HP^{-1}|S(N';h')|\ll\sum_{n\in(N'-H/P,N']}\;\max_{j\le 2H/P}|S(n;j)|.\]
Inserting this bound into (\ref{b1}) we find that
\[|S(N;h)|\ll (\log P)H^{-1}\sum_n A(n;N)\max_{j\le 2H/P}|S(n;j)|,\]
where
\begin{eqnarray*}
A(n,N):&=&\#\{(a,p):\,P<p\le 2P,\, 0\le a<p,\, n\le N'<n+H/P\},\\
&=&\#\{(a,p):\, n\le (N-aq)/p <n+H/P\}.
\end{eqnarray*}
Since
\[\sum_n A(n,N)=\sum_{a,p}\#\{n:n\le N'<n+H/P\}\le \sum_{a,p}\frac{H}{P}\ll PH\]
we deduce from H\"{o}lder's inequality that
\[|S(N;h)|^r\ll (\log P)^rP^{r-1}H^{-1}\sum_n A(n;N)\max_{j\le 2H/P}|S(n;j)|^r,\]
for any $h\le H$.  It should be noted that $A(n,N)=0$ unless $|n|\le
2q$, so that the sum over $n$ may be restricted to this range.

We proceed to sum over the values $N=N_j$ in Lemma \ref{main}, finding that
\[\sum_{j=1}^J\max_{h\le H}|S(N_j;h)|^r\ll (\log P)^rP^{r-1}H^{-1}
\sum_n A(n)\max_{j\le 2H/P}|S(n;j)|^r,\]
where
\[A(n):=\#\{(a,p,N_j):\, n\le (N_j-aq)/p<n+H/P\}.\]
From Cauchy's inequality we then deduce that
\[\sum_{j=1}^J\max_{h\le H}|S(N_j;h)|^r\ll (\log P)^rP^{r-1}H^{-1}
\mathcal{N}^{1/2}
\left\{\sum_{|n|\le 2q}\max_{j\le 2H/P}|S(n;j)|^{2r}\right\}^{1/2},\]
where
\[\mathcal{N}:=\sum_n A(n)^2\le HP^{-1}\mathcal{M},\]
with
\[\mathcal{M}:=\#\{(a_1,a_2,p_1,p_2,N_j,N_k):\,
|(N_j-a_1q)/p_1-(N_k-a_2q)/p_2|\le H/P\}.\]
Thus
\begin{eqnarray*}
\lefteqn{\sum_{j=1}^J\max_{h\le H}|S(N_j;h)|^r}\\
&\ll& (\log P)^rP^{r-3/2}H^{-1/2}
\mathcal{M}^{1/2}
\left\{\sum_{|n|\le 2q}\max_{j\le 2H/P}|S(n;j)|^{2r}\right\}^{1/2}.
\end{eqnarray*}
The second sum on the right may be bounded via Lemma \ref{l2}, giving
\[\sum_{j=1}^J\max_{h\le H}|S(N_j;h)|^r\ll_{\ep,r} 
 q^{\ep}P^{r-3/2}H^{-1/2}(q^{1/2}(H/P)^{r/2}+q^{1/4}(H/P)^r)\mathcal{M}^{1/2},\]
on replacing $\ep$ by $\ep/2$.

Naturally, in order to apply Lemma \ref{l2} we will need to have $q$
cube-free, or $r\le 3$. The
natural choice for $P$ is to take $2Hq^{-1/(2r)}\le P\ll Hq^{-1/(2r)}$ so that
$q^{1/2}(H/P)^{r/2}$ and $q^{1/4}(H/P)^r$ have the same order of
magnitude. The conditions previously imposed on $P$ are then satisfied
provided that $H\ge
q^{1/(2r)}$. With this choice for $P$ we deduce that
\begin{equation}\label{big}
\sum_{j=1}^J\max_{h\le H}|S(N_j;h)|^r\ll_{\ep,r} 
q^{1/4+3/(4r)+\ep}H^{r-2}\mathcal{M}^{1/2}.
\end{equation}

\section{Estimating $\mathcal{M}$}

In this section we will estimate $\mathcal{M}$ and complete the proof
of Lemma \ref{main}. It is the treatment of $\mathcal{M}$ which represents
the most novel part of our argument.

We split $\cM$ as $\cM_1+\cM_2$ where $\cM_1$ counts solutions with
$p_1=p_2$ and $\cM_2$ corresponds to $p_1\not=p_2$.
When $p_1=p_2$ the defining condition for $\cM$ becomes 
\[|(N_j-N_k)-q(a_1-a_2)|\le p_1H/P\le 2H.\]
Thus 
\[|a_1-a_2|\le q^{-1}(|N_j-N_k|+2H)\le 3.\]
Moreover, given $N_k$ and $a_1-a_2$, there will be at most 5 choices
for $N_j$, in view of the spacing condition (\ref{space}). Thus we
must allow for $O(P)$ choices for $p_1$, for $O(P)$ choices for $a_1$
and $a_2$, and $O(J)$ choices for $N_j$ and $N_k$, so that
\begin{equation}\label{M1}
\cM_1\ll P^2J.
\end{equation}

To handle $\cM_2$ we begin by choosing a prime $\ell$ in the range
\[q/H<\ell\le 2q/H.\]  
This is possible, by Bertrand's Postulate.  We
then set
\[M_j:=\left[\frac{N_j\ell}{q}\right],\;\;\;(1\le j\le J)\]
so that the $M_j$ are non-negative integers in $[0,\ell)$.  Moreover
the spacing condition (\ref{space}) implies that
\[M_{j+1}>\frac{N_{j+1}\ell}{q}-1\ge \frac{(N_j+H)\ell}{q}-1>
\frac{N_j\ell}{q}\ge M_j,\]
so that the integers $M_j$ form a strictly increasing sequence.
Since 
\[|N_j-qM_j/\ell|\le q/\ell\]
we now see that if
$(a_1,a_2,p_1,p_2,N_j,N_k)$ is counted by $\cM_2$ then
\[\left|\frac{qM_j/\ell-a_1q}{p_1}-\frac{qM_k/\ell-a_2q}{p_2}\right|
\le \frac{H}{P}+\frac{q}{\ell p_1}+\frac{q}{\ell p_2},\]
whence
\[|p_2M_j-p_1M_k-\ell\delta|\le\frac{H\ell p_1p_2}{Pq}+p_1+p_2\le
12P,\]
with $\delta=a_1p_2-a_2p_1$.
If $p_1,p_2$ and $\delta$ are given, there is at most one pair of 
integers $a_1,a_2$ with $0\le a_1<p_1$, $0\le a_2<p_2$ and 
$a_1p_2-a_2p_1=\delta$. Thus
\[\cM_2\le\sum_{M_j,M_k}\#\{(p_1,p_2,m):\, |m|\le 12P, 
\, p_2M_j-p_1M_k\equiv m\mod{\ell}\}.\]

We now consider how many pairs $p_1,p_2$ there may be for each choice
of $M_j,M_k$. We define the set
\[\Lambda:=\{(x,y,z)\in\Z^3:\, xM_j-yM_k\equiv z\mod{\ell}\},\]
which will be an integer lattice of determinant $\ell$. Admissible
pairs $p_1,p_2$ produce points $\x=(x,y,z)\in\Lambda$ with
$x\not=y$ both prime and $|\x|\le 12P$, where 
\[|\x|:=\max(|x|,|y|,|z|). \]
The lattice $\Lambda$ has a $\Z$-basis
$\b{b}_1,\b{b}_2,\b{b}_3$ such that 
\begin{equation}\label{lat1}
|\b{b}_1|\le |\b{b}_2|\le |\b{b}_3|
\end{equation}
and
\begin{equation}\label{lat2}
\det(\Lambda)\ll |\b{b}_1|.|\b{b}_2|.|\b{b}_3|\ll\det(\Lambda)=\ell,
\end{equation}
and with the property that there is an absolute constant $c_0$ such that
if $\x\in\Lambda$ is written as $\lambda_1\b{b}_1+\lambda_2\b{b}_2+
\lambda_3\b{b}_3$ then  
\[|\lambda_i|\le c_0|\x|/|\b{b}_i|,\;\;\; (1\le i\le 3).\]
The existence of such a basis is a standard fact about lattices, see
Browning and Heath-Brown \cite[Lemma 1, (ii)]{bhb}, for example. When
$|\b{b}_3|\le 12c_0P$ we
now see that the number of lattice elements of size at most $12P$ is
\begin{eqnarray*}
&\le&\left(1+\frac{12c_0P}{|\b{b}_1|}\right)
\left(1+\frac{12c_0P}{|\b{b}_2|}\right)
\left(1+\frac{12c_0P}{|\b{b}_3|}\right)\\
&\ll &\frac{P^3}{|\b{b}_1|.|\b{b}_2|.|\b{b}_3|}\\
&\ll &\frac{P^3}{\det(\Lambda)}\\
&\ll & HP^3q^{-1}
\end{eqnarray*}
by (\ref{lat1}) and (\ref{lat2}). If $|\b{b}_1|>12c_0P$ the only
vector in $\Lambda$ of norm at most $12P$ is the origin, while if
$|\b{b}_1|\le 12c_0P<|\b{b}_2|$ the only possible vectors are of the
form $\lambda_1\b{b}_1$. In this latter case
$(p_2,p_1,m)=\lambda_1\b{b}_1$ so that $\lambda_1$ divides ${\rm
  h.c.f.}(p_2,p_1)=1$. Hence there is at most 1 solution in this
case. 

There remains the situation in which $|\b{b}_2|\le 12c_0P<|\b{b}_3|$,
so that the admissible vectors are linear combinations
$\lambda_1\b{b}_1+\lambda_2\b{b}_2$.  In this case we write 
$\b{b}_i=(x_i,y_i,z_i)$ for $i=1,2$ and set $\Delta=x_1y_2-x_2y_1$.
If $\Delta=0$ then $(x_1,y_1)$ and $(x_2,y_2)$ are proportional, and
hence are both integral scalar multiples of some primitive vector
$(x,y)$ say.  However we then see that if
$(p_2,p_1,m)=\lambda_1\b{b}_1+\lambda_2\b{b}_2$ then $(p_2,p_1)$ is a
scalar multiple of $(x,y)$, so that $\b{b}_1$ and $\b{b}_2$ determine
$p_1$ and $p_2$.  Thus when $\Delta=0$ the primes $p_1$ and $p_2$ are
determined by $M_j$ and $M_k$. In order to summarize our conclusions
up to this
point we write $\cM_3$ for the contribution to $\cM_2$ corresponding to all
cases except that in which $|\b{b}_2|\le 12c_0P<|\b{b}_3|$ and
$\Delta\not=0$. With this notation we then have
\begin{equation}\label{M2}
\cM_3\ll (HP^3q^{-1}+1)J^2.
\end{equation}

Suppose now that $|\b{b}_2|\le 12c_0P<|\b{b}_3|$ and $\Delta\not=0$.
We will write $\cM_4$ for the corresponding contribution to $\cM$.
In this case we must have $\lambda_3=0$, and the number of choices for 
$\lambda_1$ and $\lambda_2$ will be
\[\le\left(1+\frac{12c_0P}{|\b{b}_1|}\right)
\left(1+\frac{12c_0P}{|\b{b}_2|}\right)
\ll \frac{P^2}{|\b{b}_1|.|\b{b}_2|}.\]
Thus if $L<|\b{b}_1|.|\b{b}_2|\le 2L$, say, the contribution to $\cM_4$
will be $O(P^2L^{-1})$ for each pair $M_j,M_k$.  

To estimate the number of pairs of vectors $\b{b}_1,\b{b}_2$ with
$L<|\b{b}_1|.|\b{b}_2|\le 2L$ we observe that there are
$O(B_1^3B_2^3)$ possible choices with
$B_1<|\b{b}_1|\le 2B_1$ and $B_2<|\b{b}_2|\le 2B_2$.  A dyadic
subdivision then shows that we will have to consider $O(L^3\log L)$ pairs
$\b{b}_1,\b{b}_2$.  Writing $\b{b}_i=(x_i,y_i,z_i)$ for $i=1,2$ as
before we will have
\[x_1M_j-y_1M_k\equiv z_1\mod{\ell},\;\;\;
x_2M_j-y_2M_k\equiv z_2\mod{\ell}.\]
These congruences determine $\Delta M_j$ and $\Delta M_k$ modulo
$\ell$, and since $\ell$ is prime and $0\le M_j,M_k<\ell$ we see
that $\b{b}_1$ and $\b{b}_2$ determine $M_j,M_k$ precisely, providing
that $\ell\nmid\Delta$. However 
\[|\Delta|\le 2|\b{b}_1||\b{b}_2|\le
2(|\b{b}_1|.|\b{b}_2|.|\b{b}_3|)^{2/3}\ll\det(\Lambda)^{2/3}=\ell^{2/3}\]
by (\ref{lat1}) and (\ref{lat2}).  Since $\Delta\not=0$ we then see
that $\ell\nmid\Delta$ providing that $q/H$, or equivalently
$\ell$, is sufficiently large.  Under this assumption we therefore
conclude that there are $O(L^3\log L)$ pairs $M_j,M_k$ for which
$|\b{b}_2|\le 12c_0P<|\b{b}_3|$ and $\Delta\not=0$ and for which
$L<|\b{b}_1|.|\b{b}_2|\le 2L$.  Thus each dyadic range $(L,2L]$ contributes
$O(P^2L^{-1}\min(J^2,L^3\log L))$ to $\cM_4$. Since
\[P^2L^{-1}\min(J^2,L^3)\le
P^2L^{-1}(J^2)^{2/3}(L^3)^{1/3}=P^2J^{4/3}\]
we deduce that
\[\cM_4\ll P^2J^{2/3}\log q,\]
and comparing this with the bounds (\ref{M1}) and (\ref{M2}) we then
see that
\[\cM\ll (HP^3q^{-1}+1)J^2+P^2J^{4/3}\log q.\]

We may now insert this bound into (\ref{big}), recalling that $P$ is
of order $Hq^{-1/(2r)}$ to deduce, after replacing $\ep$ by $\ep/2$ that
\begin{eqnarray*}
\lefteqn{\sum_{j=1}^J\max_{h\le H}|S(N_j;h)|^r}\\
&&\ll_{\ep,r} 
q^{1/4+1/(4r)+\ep}H^{r-1}\{J^{2/3}+J(H^{-1}q^{1/(2r)}+Hq^{-1/2-1/(4r)})\},
\end{eqnarray*}
as required for Lemma \ref{main}.

\section{Deduction of the theorem}

We will prove the theorem by induction on $r$. The result for $r=1$ is
an immediate consequence of Lemma \ref{l3}, together with the
P\'{o}lya--Vinogradov inequality.

For $r\ge 2$ we will use a dyadic subdivision, classifying the $N_j$
according to the value $V=2^v$ for which 
\begin{equation}\label{dy}
V/2<\max_{h\le H}|S(N_j;h)|^r\le V.
\end{equation}
Clearly numbers $N_j$ for which the corresponding $V$ is less than $1$
make a satisfactory contribution in our theorem, and so it suffices to
assume that (\ref{dy}) holds for all $N_j$.

We now give three separate arguments, depending on which of the three
terms on the right of (\ref{est}) dominates.  If
\[\sum_{j=1}^J\max_{h\le H}|S(N_j;h)|^r\ll_{\ep,r} 
q^{1/4+1/(4r)+\ep}H^{r-1}J^{2/3}\]
then
\[JV^r\ll_{\ep,r} q^{1/4+1/(4r)+\ep}H^{r-1}J^{2/3},\]
whence
\[JV^{3r}\ll_{\ep,r} q^{3/4+3/(4r)+3\ep}H^{3r-3},\]
which suffices for the theorem.  If the second term dominates we will
have 
\[\sum_{j=1}^J\max_{h\le H}|S(N_j;h)|^r\ll_{\ep,r} 
q^{1/4+1/(4r)+\ep}H^{r-1}JH^{-1}q^{1/(2r)},\]
so that
\[JV^r\ll_{\ep,r} q^{1/4+3/(4r)+\ep}H^{r-2}J.\]
In this case it follows that
\begin{equation}\label{by}
V^r\ll_{\ep,r} q^{1/4+3/(4r)+\ep}H^{r-2}.
\end{equation}
We now use Lemma \ref{l3}, which implies that
\begin{equation}\label{bz}
JV^{2r}\ll_{\ep,r} q^{\ep}(qH^{r-1}+q^{1/2}H^{2r-1})
\ll_{\ep,r} q^{1/2+\ep}H^{2r-1}
\end{equation}
since $H\ge q^{1/(2r)}$. Coupled with (\ref{by}) this yields
\[JV^{3r}\ll_{\ep,r} q^{3/4+3/(4r)+2\ep}H^{3r-3}\]
which again suffices for the theorem.
Finally, if the third term on the right of (\ref{est}) dominates we
must have 
\[JV^r\ll_{\ep,r} q^{-1/4+\ep}H^rJ\]
whence $V\ll_{\ep}Hq^{-1/(4r)+\ep/r}$. Here we shall
use the inductive hypothesis, which tells us that
\[JV^{3r-3}\ll_{\ep,r} q^{3/4+3/(4r-4)+\ep}H^{3r-6}\]
if either $r=2$ or $H\ge q^{1/(2r-2)}$ and $r\ge 3$.  Under this
latter assumption we therefore deduce that
\[JV^{3r}\ll_{\ep,r}q^{3/4+\phi+4\ep}H^{3r-3}\]
with
\[\phi=\frac{3}{4r-4}-\frac{3}{4r}\le \frac{3}{4r}\]
for $r\ge 2$.  It therefore remains to consider the case in which
$r\ge 3$ and $q^{1/(2r)}\le H\le q^{1/(2r-2)}$. However for such $H$
we may again use the bound (\ref{bz}), whence
\begin{eqnarray*}
JV^{3r}&\ll_{\ep,r}&q^{1/2+\ep}H^{2r-1}.\,q^{-1/4+\ep}H^r\\
&=& q^{1/4+2\ep}H^{3r-3}\{Hq^{-1/(2r-2)}\}^2q^{1/(r-1)}\\
&\le &q^{1/4+1/(r-1)+2\ep}H^{3r-3}.
\end{eqnarray*}
To complete the proof of this final case it remains to observe that
$1/4+1/(r-1)\le 3/4+3/(4r)$ for $r\ge 3$.

\bigskip
\bigskip

Mathematical Institute,

24--29, St. Giles',

Oxford

OX1 3LB

UK
\bigskip

{\tt rhb@maths.ox.ac.uk}

\end{document}